\documentclass[preprint, 1p, 12pt]{elsarticle}

\usepackage{amsmath,amsthm,amssymb}
\usepackage{graphicx}
 \usepackage{algorithm} 
 \usepackage{longtable} 
 \usepackage[noend]{algpseudocode} 
\makeatletter
\def\BState{\State\hskip-\ALG@thistlm}
\makeatother
\usepackage{mathrsfs}
 \usepackage{color}
\newcommand{\nmu}{\hat{\mu}}
\usepackage{booktabs} 
\makeatletter
\newcommand\footnoteref[1]{\protected@xdef\@thefnmark{\ref{#1}}\@footnotemark}
\makeatother

\usepackage{tikz}
\usetikzlibrary{shadows}
\usetikzlibrary{intersections}
\tikzset{
	shadowed/.style={preaction={transform canvas={shift={(1pt,-1pt)}},draw=gray,very thick, opacity=0.3}},
  }

\definecolor{darkblue}{HTML}{AAC4DC}
\definecolor{maniblue}{HTML}{A7BAD1}
\definecolor{myred}{HTML}{B18D94}
\definecolor{myblue}{HTML}{4F81BD}
\begin{document}
\begin{frontmatter}
\title{Alternatives for Generating a Reduced Basis to Solve the Hyperspectral Diffuse Optical Tomography Model}
\author[gc]{Rachel~Grotheer\corref{cor1}}
\ead{rachel.grotheer@goucher.edu}
\author[b]{Thilo~Strauss}
\ead{thilum@gmx.de} 
\author[ub]{Phil~Gralla}
\ead{phil.gralla@uni-bremen.de}
 \author[cu]{Taufiquar~Khan}
 \ead{khan@clemson.edu}  
\cortext[cor1]{Corresponding author}

\address[gc]{Goucher College, Center for Data, Mathematical, and Computational Sciences, Baltimore, MD, 21204, USA } 
\address[b]{University of Washington, Department of Pediatrics, 1959 NE Pacific St, Seattle, WA 98195, USA}
\address[ub]{University of Bremen, Zentrum f{\"u}r Technomathematik (ZeTeM), Bibliothekstra{\ss}e 5, 28359 Bremen, Germany}
\address[cu]{Clemson University, Department of Mathematical Sciences, Clemson, SC, 29634 USA}

%

\begin{abstract} The Reduced Basis Method (RBM)\footnote{\label{fn1}Abbreviations used in this article: Reduced Basis Method (RBM), hyperpsectral diffuse optical tomography (hyDOT)} is a model reduction technique used to solve parametric PDEs that relies upon a basis set of solutions to the PDE at specific parameter values.
To generate this reduced basis, the set of a small number of parameter values must be strategically chosen. 
 We apply a Metropolis algorithm and a gradient algorithm to find the set of parameters and compare them to the standard greedy algorithm most commonly used in the RBM. 
We test our methods by using the RBM to solve a simplified version of the governing partial differential equation for hyperspectral diffuse optical tomography (hyDOT)\footnoteref{fn1}. 
The governing equation for hyDOT is an elliptic PDE parameterized by the wavelength of the laser source.
For this one-dimensional problem, we find that both the Metropolis and gradient algorithms are potentially superior alternatives to the greedy algorithm in that they generate a reduced basis which produces solutions with a smaller relative error with respect to solutions found using the finite element method and in less time. 
\end{abstract}

 \end{frontmatter}
\section{Introduction}
Since its development in the late 1970s~\cite{almroth,nagy,noor}, the reduced basis method (RBM) has been used as a model reduction technique to approximate the solution to a large dimensional problem in a variety of applications. 
For example, it has been applied to the nonlinear structural analysis of beams and arches~\cite{almroth, nagy, noor},  to differential algebraic equation systems~\cite{porsching87}, and to solving a variety of parametric PDEs,  including the linear Helmholtz-elasticity equation~\cite{cuong, rozza} and the nonlinear Navier-Stokes equations~\cite{gunzburgerbook, peterson83, peterson, porsching85}.
It has also been applied in solving eigenvalue problems~\cite{machiels,horger16} and control problems~\cite{ito98,ito98reduced, negri13, tonn11}. 
Most recently, the RBM has had success in accurately approximating the solution of parameterized elliptic partial differential equations (PDEs), especially linear elliptic PDEs with an affine parameter dependence~\cite{maday,veroy,veroy02}. 
The RBM is particularly effective in providing efficient, real-time approximations to the solution of parametric PDEs, presented in their weak form as an input-output functional
\begin{equation}
a(u(\mu), v; \mu) = f(v).
\label{eq:inputoutput}
\end{equation}

Since the RBM only requires the ``truth" solution for select parameter values, it offers significant computational reduction, 
an advantage over a model reduction technique like proper orthogonal decomposition, for example. 

In recent years, optical imaging, imaging using a low-energy light source in the visible to near infrared range, has arisen as a potential alternative to traditional medical imaging techniques such as x-ray or computed tomography (CT). 
The advantages of optical imaging are that the low-energy light is non-ionizing and thus, not harmful to tissue, the devices cost less than existing medical imaging devices, and they are helpful in providing functional, rather than anatomical, information~\cite{arridgereview}. 
A popular form of optical imaging is known as Diffuse Optical Tomography (DOT). 

Currently, researchers are seeking to apply hyperspectral imaging, where hundreds of optical wavelengths are used for the laser source, to DOT to create what is known as hyperspectral DOT, or hyDOT~\cite{larusson,larusson_pals,lufei,saibaba,saibaba15}. 
When hyperspectral imaging is applied to DOT, the governing PDE is parameterized with respect to the wavelength, $\lambda$. 
This parameterization increases the dimension of the problem and makes it an ideal candidate for use of model reduction techniques in order to approximate the solution. 

Potential model reduction techniques include Proper Orthogonal Decomposition~\cite{dugunzburger,sirovich} and the Empirical Interpolation Method~\cite{barrault, lassila}. 
The Empirical Interpolation Method is often used in conjunction with the RBM, especially to serve as an approximate replacement to the affine decomposition assumption (discussed below) when it cannot be met~\cite{casenave14}. 
While there are other model reduction techniques and the effectiveness of these methods is well-established, we posit that the reduced basis method will be especially effective when applied to the hyDOT image reconstruction problem due to its ill-posed nature.  
We note that only numerical results of the effectiveness of the applying the RBM to the forward problem of hyDOT are presented in this paper. A more rigorous theoretical explanation and further numerical results are introduced in~\cite{dissertation}. 

Application of the RBM to a parameterized elliptic PDE, such as the governing equation for hyDOT, requires sampling the parameter space in order to construct a reduced basis with which to obtain an approximation of the solution. 
If the basis generating set is chosen wisely, it has been shown that the approximation converges to the true solution exponentially with respect to the number, $N$, of basis functions in the case where the parameter defining the parameterization is one-dimensional~\cite{buffa,maday02,maday}, and has been experimentally verified in higher dimensions~\cite{rozza,sen08,veroy}. 
As a result, there has been much work done in recent years to develop efficient and effective algorithms for choosing this set of parameters. 
Most of the recent approaches use a greedy algorithm developed using \textit{a posteriori} error bounds~\cite{buffa,hesthaven14,sen08,veroy,veroy02}. 
Alternatively, in recent years there has been a vast amount of research in the theory and application of different Metropolis Hastings Algorithms \cite{bardsley2012mcmc, Adlouni2006, Gelman1996, Raftery1995, thilo}.
Because of the relatively slow convergence rate, researchers created adaptive Metropolis Algorithms \cite{Gilks1998, Haario2001}, which are considerably faster than the original version developed in 1953 by Metropolis et al. \cite{metropolis1953equation} and later generalized in 1970 by Hastings et al. \cite{hastings1970monte}. Here, we propose a special kind of adaptive Metropolis Algorithm to estimate parameters to generate a reduced basis and compare its effectiveness and speed to the traditional greedy algorithm and to a gradient algorithm.
We test the algorithms on an application of the RBM to the forward problem in hyDOT.
The governing PDE is the first order diffusion approximation of the radiative transport equation. 
The one-dimensional parameter in question is the optical wavelength, $\lambda$, of the laser source.

The article is organized as follows. Section 2 explains the Reduced Basis Method and describes the standard greedy algorithm used in most applications of the RBM, a gradient algorithm, and the pilot adaptive Metropolis Hastings Algorithm introduced in~\cite{thilo}, and their application to finding the basis generating set of parameters for the RBM. 
Then, Section 3 gives an overview of hyperspectral diffuse optical tomography. Section 4 gives a description of the simplified one-dimensional hyDOT simulation (where the PDE is parameterized by the wavelength, $\lambda$) on which all three algorithms were tested for a variety of basis sizes. A comparison of the total error in the solution and the computational time between the three algorithms is given. 
Section 5 gives a summary of the results and directions for future work. 

\section{Reduced Basis Method}
\begin{figure}[h]
\centering
\begin{tikzpicture}
 \draw[myblue,shadowed,thick,->] (0,0) -- (0,3);
 \draw[myblue,rotate = -120,shadowed, thick, ->] (0,0) -- (0,3);
 \draw[myblue,rotate = 120, shadowed, thick,->] (0,0) -- (0,3);
 \draw[thick,shadowed] plot [smooth, tension=1] coordinates {(-2.5,0) (-1,1.5) (.7,.5) (2.3,2.5)};
 \node at (-2,1.7) {$\mathcal{M}$};
 \node at (0,-1) {X($\Omega$)};
\end{tikzpicture}
\begin{tikzpicture}
 \fill[darkblue!40!white] (-2.5,.2)--(-1.5,2.7)--(3,2.7)--(2,.2)--(-2.5,.2);
 \draw[myblue,shadowed,thick,->] (0,0) -- (0,3);
 \draw[myblue,rotate = -120,shadowed, thick, ->] (0,0) -- (0,3);
 \draw[myblue,rotate = 120, shadowed, thick,->] (0,0) -- (0,3);
 \draw[thick,shadowed, name path = graph2] plot [smooth, tension=1] coordinates {(-2.5,0) (-1,1.5) (.7,.5) (2.3,2.5)};
 \draw[red, shadowed, name path= graph1] plot [smooth, tension=1] coordinates {(-2.5,.2) (-1,1.3) (.7,.6) (2.3,2)};
 \node at (-2,1.7) {$\mathcal{M}$};
 \node at (0,-1) {X($\Omega$)};
 \node at (2.3,1.3){$W_n$};
 \fill [myblue, name intersections={of=graph1 and graph2, total=\n}]  
      \foreach \i in {2,...,\n}{(intersection-\i) circle [radius=2pt]};
 \draw [name intersections={of=graph1 and graph2, total=\n}]  
      \foreach \i in {1,...,\n} {(intersection-\i) coordinate (my-intersection-\i)};
 \node[below right] at (my-intersection-2) {$u(\mu_1)$};
 \node[below right] at (my-intersection-4) {$u(\mu_n)$};
 \path [name path = X] (-1,0)--(-1,3);
 \path [name intersections={of=graph1 and X,by=S}];
 \node [label=90:\footnotesize{$u(\mu)$}] at (S) {};
 \draw[thick, myblue] (S) circle (.08); 
\end{tikzpicture}
\caption[Diagram representing the idea of the manifold in the RBM]{The graph on the left represents the finite element space, $X(\Omega)$ of a PDE with solution $u$ and domain $\Omega$, and  a lower-dimensional manifold, $\mathcal{M}$, on which the reduced basis parameter-induced solutions reside. On the right, the solution $u$ at unknown parameter $\mu$ is approximated by the RBM using a linear combination of basis functions for $\mathcal{M}$, $\{u(\mu_i)\}_{i=1}^N$ for $N \ll \text{dim}(X)$. This approximated solution lies on the subspace $W_N$ approximating the manifold, shown in red. }
\label{fig:manifoldnew}
\end{figure}
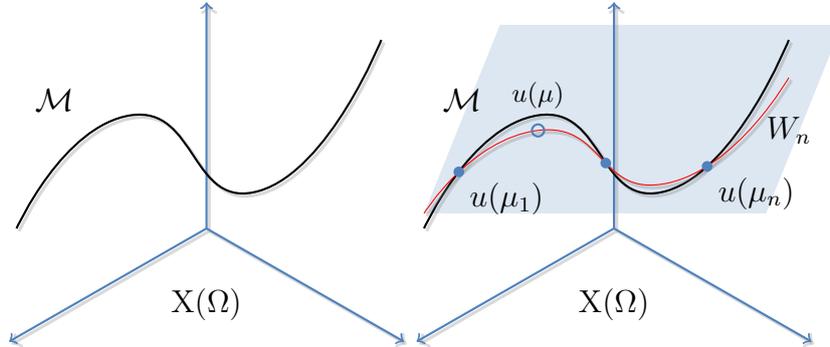

The key idea of reduced basis methods is to take a high-dimensional, computationally expensive problem, and reduce it to a problem on a space of lower dimension that is determined by a select number of solutions to the high-dimensional problem. 
These methods are well-suited to parametric partial differential equations since the parametric-dependence of the state variable means that it actually resides (or ``evolves") on a parameter-induced manifold that is of much lower dimension than the infinite dimensional space associated with the PDE. 
Thus, the RBM seeks to approximate the solution, $u(\mu)$, to a PDE at any parameter, $\mu$, by creating a representative subspace, $W_N$, approximating the manifold~\cite{lassila_bookchap}. This subspace is created by finding the span of a small collection of high-dimensional solutions at a selected parameter values $\mu_1, \ldots, \mu_N$, which are the basis elements of this subspace. 
An illustration of this is given in Figure~\ref{fig:manifoldnew}. 

In general, the solution $u(x;\mu)$ to a parameter-dependent PDE with parameter $\mu$, 
is discretized and approximated to give what is known as a ``truth approximation" of the solution
\begin{equation}
u^{fe}(\mu) = \sum_{i=1}^\mathcal{N} c_i(\mu)\phi_i 
\label{eq:truthapprox}
\end{equation}
where the $\phi_i$ form a basis for the finite element approximation space, $X$, of dimension $\mathcal{N}$, where $\mathcal{N}$ is large~\cite{maday}. 
Note that $X$ is a Hilbert space over the bounded spatial domain $\Omega \subseteq \mathbb{R}^d$. 
Thus, the reduced basis approximation is not an approximation to the infinite dimensional problem, but rather an approximation to its finite element discretization, which is still of large and unmanageable dimension~\cite{rozza}. 
From its weak formulation, the PDE is defined as an input-output operator given as a bilinear form. 
In general, for input parameter $\mu$ (where $\mu$ may be a vector if there is more than one parameter), this is given as
\begin{equation}
a(u(\mu), v; \mu) = f(v), \qquad \forall v \in X 
\label{eq:genPDE}
\end{equation}
The bilinear (or sesquilinear in the complex case) form is continuous with respect to $X$, as well as coercive. 
Often, it is also assumed to be symmetric in the real case to prove convergence of the reduced basis approximate solution to the finite element solution~\cite{cuong_thesis}. 
The output is a linear functional of the state variable $u(\mu)$ and is given by $s(\mu) = \ell(u(\mu))$. 
Problems where $\ell = f$ are commonly called ``compliant"; otherwise, the problem is known as ``noncompliant." 
The set of all possible solutions to the parameterized PDE is a manifold of dimension $N$ given by
\[\mathcal{M} = \{ u(\mu) | \mu \in \mathcal{D} \subset \mathbb{R} \} \]
where $\mathcal{D}$ is the parameter space. 

To construct a basis for $\mathcal{M}$, first a set of parameters $S_N = \{ \mu_i| i = 1, \ldots, N\}$ is chosen, typically using a greedy algorithm like the one described in Section 2, choosing the parameters that give the most significant information about the structure of the manifold. 
Note that $N \ll \mathcal{N}$.
From this set of parameters, we define a basis space for the manifold,
\[ W_N = \text{span} \{ u(\mu_i) | i \leq N\}. \]
where the $u(\mu_i)$ are defined by the finite element approximation given in (\ref{eq:truthapprox}). 
Once the approximation space is constructed, the reduced basis approximation of the solution for a new $\mu$ is given by
\begin{equation}
u_N(\mu) = \sum_{j=1}^N \hat{c}_j(\mu)u(\mu_j) \in W_N 
\label{eq:PDEapprox}
\end{equation}
The coefficients $\hat{c}_j$ are found by using a Galerkin or Petrov-Galerkin projection,
 with respect to $W_N$ to find a solution to 
  \begin{eqnarray*}
  a\left(\sum_{j=1}^N \hat{c}_j(\mu)u(\mu_j),v; \mu\right) &=& f(v), \qquad \forall v \in W_N\\
  \implies \sum_{j=1}^N \hat{c}_j(\mu) a(u(\mu_j), v; \mu) &=& f(v).
  \end{eqnarray*}
  Since $W_N = \text{span}\{u(\mu_i)| i \leq N\}$, $a(\cdot, \cdot)$ is bilinear and $f(\cdot)$ is linear, it is sufficient to solve 
  \begin{equation}
   \sum_{i=1}^N \hat{c}_j(\mu)a(u(\mu_i), u(\mu_j); \mu) = f(u(\mu_j)), \qquad \forall j \in \{1, 2, \ldots, N\}.
   \label{eq:galerkin} 
   \end{equation}
This problem is ill-conditioned because the basis vectors $u(\mu_i)$ are usually pointing in similar directions due to the smoothness of the low-dimensional $\mathcal{M}$ and thus are not sufficient to solve for the unknown $\hat{c}_i$. 
Therefore, the Gram-Schmidt procedure is often applied to make the $u(\mu_i)$ orthogonal.

The RBM has both an offline stage and an online stage when it is implemented computationally. 
In order for the RBM to have the full benefit of computational efficiency, as much as possible should be precomputed and stored in the offline stage. 
In this stage, the sample set $S_N$ is constructed, and the basis functions $u(\mu_i)$ of the space $W_N$ are precomputed using the FEM approximation.

A significant source of computational effort for the RBM is the dependence of the bilinear form $a(u(\mu_i), u(\mu_j); \mu)$ on the parameter $\mu$, which results in it having to be calculated for each pair $(u(\mu_i), u(\mu_j))$ in the Galerkin projection (\ref{eq:galerkin}) during the online stage. 
If that parameter dependence can be removed, then the bilinear form can be precomputed and stored in the offline stage. 
For some PDEs, the bilinear form can be decomposed as
\begin{equation}
a(u, v; \mu) = \sum_{q=1}^{Q _a}\Theta^q(\mu)a^q(u,v) 
\label{eq:affinedecomp}
\end{equation}
where the $\Theta^q: \mathcal{D} \to \mathbb{R}$ are differentiable and, in general, very smooth functions depending on $\mu$, and the $a^q: X \times X \to \mathbb{R}$ are parameter-independent, continuous (with respect to $X$) bilinear forms~\cite{rozza,cuong_thesis}. 
Given this decomposition, the reduced basis approximation involves solving (plugging~(\ref{eq:affinedecomp}) into~(\ref{eq:galerkin})),
\begin{equation}
\sum_{i=1}^N \sum_{q=1}^{Q_a} \hat{c}_i(\mu)\Theta^q(\mu)a^q(u(\mu_i), u(\mu_j)) = f(u(\mu_j)) \qquad \forall j = 1, 2, \ldots, N
\label{eq:RBapproxdecomp}
\end{equation}
for the coefficients $\hat{c}_i(\mu)$ for each new value of the parameter $\mu$. 
Since they no longer depend on the parameter, the $a^q(u(\mu_i),u(\mu_j))$ can be precomputed, reducing the computational burden significantly. 
Note that $f$ may also be affinely parameter dependent, and so there may also exist a similar decomposition for $f$~\cite{rozza,lassila_bookchap}. 
For an example of how the RBM can be extended to PDE's with nonaffine parameter dependence, see~\cite{grepl}. 
The great benefit of the RBM is that all operation counts in the online stage are independent of $\mathcal{N} \gg N$ and so the computational burden is greatly reduced from a procedure that would require finding the finite element solution for each new value of $\mu$. 
Note that significant computational reduction is only possible if the affine decomposition of the bilinear form given in~(\ref{eq:affinedecomp}) is possible. 

In summary, Algorithm~\ref{alg:rbm} describes how to solve a parameterized partial differential equation using the reduced basis method:
\begin{algorithm}[H]
\caption{Reduced Basis Method}\label{rbm}
\begin{algorithmic}[1]
\Procedure{Offline Stage}{}
\BState Finite element discretization: $u^{fe} = \sum_{i=1}^\mathcal{N} c_i \phi_i$
\BState Choose parameter samples : $\mu_1, \ldots, \mu_N, N \ll \mathcal{N}$
\BState Define $W_N = \text{span}\{u(\mu_k), k = 1, \ldots, N\}$, where $u(\mu_k) = \sum_{i=1}^\mathcal{N} c_i(\mu_k)\phi_i$
\BState Compute and store $\sum_{i=1}^N a^q(u(\mu_i), u(\mu_j))$, for $q = 1, \ldots, Q^a$, $j = 1, 2, \ldots, N$ 
\EndProcedure
\Procedure{Online Stage}{}
\BState Find a solution with respect to $\hat{c}_i$ for \newline
$\sum_{i=1}^N \sum_{q=1}^Q \hat{c}_i(\mu)\Theta^q(\mu)a^q(u(\mu_i), u(\mu_j)) = f(u(\mu_j)) \qquad \forall j = 1, 2, \ldots, N$
\BState Reduced Basis approximation: $u_N(\mu) = \sum_{i=1}^N \hat{c}_i u(\mu_i)$
\EndProcedure
\end{algorithmic}
\label{alg:rbm}
\end{algorithm}
%
\subsection{Greedy Algorithm}
A key element in the RBM is choosing an appropriate sample set $S_N$ of parameters from which to form a basis for the approximation space $W_N$. 
The goal is to choose the parameters that yield the most sensitive solutions, that is, solutions that illustrate the most significant features of the solution space.
Enough parameter samples must be chosen to yield reduced basis approximations that converge to the truth solution (which in turn converges to the exact solution), but the number of samples must be far fewer than the dimension of the finite element approximation space and small enough to preserve computational efficiency. 

The most common sampling procedure for generating a reduced basis space is a greedy algorithm based on \textit{a posteriori} error bounds that allows for an efficient sampling of $\mathcal{M}$ that is independent of $\mathcal{N}$~\cite{cuong, rozza,cuong_thesis, sen08}.
First, the smallest anticipated error tolerance, $\epsilon_{\text{tol,min}}$ 
is calculated \textit{a priori} offline 
where $\epsilon_{\text{tol, min}} = \min \epsilon_\text{tol}$ and 
\[ || u^{fe}(\mu) - u_N(\mu) ||_X \leq \epsilon_{\text{tol}}. \]
Next, a very fine mesh, $\Xi$, of the parameter space $\mathcal{D}$ is created, containing the surrogate values of $\mathcal{D}$ from which the greedy algorithm will draw to generate a set of ``training" samples, $S_N$. 
This mesh is usually generated using a Monte Carlo method with respect to a uniform or log-uniform density and must be sufficiently fine to ensure that further refinement does not significantly improve the results~\cite{rozza}. 

Next, $N_\text{max}$, the maximum allowable dimension of the reduced basis space, is defined, such that the desired accuracy of the reduced basis approximation is attained~\cite{rozza,sen08}. 
$N_\text{max}$ can also be adaptively determined during the greedy algorithm~\cite{cuong}. 
A sample $\mu_1$ is chosen at random to be the first sample added to $S_1$ 
such that $S_1 = \{\mu_1\}$, and then $W_1 = \{ u(\mu_1)\}$ is calculated. 
The next sample, $\mu_2$ is calculated as
\[ \mu_2 = \text{arg max}_{\mu \in \Xi} \epsilon_1^\ast(\mu) \]
where 
\begin{equation}
\epsilon_k^\ast = \Delta_k^s (\mu)/s_k(\mu) 
\label{eq:relerror}
\end{equation}
 is the relative error bound, and $\Delta_k^s (\mu)$ is described below~\cite{cuong_thesis, sen08}. 
That is, $\epsilon_k^\ast$ is a bound on the error in approximating $u(\mu)$ by a linear combination of elements in the set $\{u(\mu_1), \ldots, u(\mu_k)\}$, $1\leq k \leq N_{\max}$. 
This approximation is denoted $u_k(\mu)$. 
Note that $s_k(\mu)$ is the reduced basis approximation of the desired output at $\mu$. 
 
 To calculate $\Delta_k^s(\mu)$, a positive lower bound, $\hat{\alpha}(\mu)$, for the inf-sup stability constant, $\alpha(\mu)$, of the bilinear form $a(u,v;\mu)$ must be determined, where 
 \[\alpha(\mu) = \inf_{w\in X} \sup_{v\in X} \frac{a(w,v; \mu)}{||w||_X||v||_X}\] 
 Then, the residual associated with $u_k(\mu)$, 
 \[ r(v;\mu) := f(v) - a(u_k(\mu), v; \mu) \qquad v \in X\]
is calculated to define the dual norm of the residual
\begin{equation}
\varepsilon_k(\mu) = ||r(\cdot ; \mu)||_{X'} = \sup_{v \in X} \frac{r(v; \mu)}{||v||_X} 
\label{eq:dualnorm}
\end{equation}
Thus, for the compliant problem, 
we define
\[ | s(\mu) - s_k(\mu) | \leq \Delta_k^s(\mu) : = \frac{\varepsilon_k^2(\mu)}{\hat{\alpha}}. \] 
A sharper bound with more rapid convergence of the reduced basis output approximation can be obtained using the dual problem for noncompliant ($\ell \neq f$) problems~\cite{sen08}. 
Thus, in calculating each new $\mu_i$, we are looking for the parameter value in $\Xi$ that will give us the largest scaled residual, that is, the parameter that affects the most significant change in the solution. 
It is important to note that $\Delta_k^s$ is reliable (an upper bound of the true error), and as a surrogate, truly represents the true error, $||u(\mu) - u_k(\mu)||_X$ (and thus is a sharp bound)~\cite{sen08}. 
 Note that other types of error estimators can be used in the definition of $\mu_i$, see, e.g. ~\cite{hesthaven14, rozza}. 

Next, it must be verified that $\varepsilon_2^\ast = \epsilon_1(\mu_2)$ is greater than $\epsilon_{\text{tol, min}}$. 
 If it is not, then the current set $S_1$ is the final sample set and the process is terminated. 
 Otherwise, $S_2 = S_1 \cup \{\mu_2\}$ and $W_2 = \text{span} \{W_1, \{ u(\mu_2) \}\}$. 
 The process then continues until $k = N_{\max}$ or an error estimate goes below $\epsilon_{\text{tol, min}}$. 
 Hesthaven, et al. have adapted this algorithm to cases where the parameter space is of high dimension~\cite{hesthaven14}. 
 The greedy algorithm described here is summarized in Algorithm~\ref{alg:ga}.
 \begin{algorithm}[h]
\caption{Greedy Sampling Algorithm}\label{alg:ga} 
\begin{algorithmic}
\For{$k = 2:N_\text{max}$}  
	\State $\mu_k = \text{arg max}_{\mu \in \Xi} \epsilon_{k-1}^\ast (\mu)$
	\State $\varepsilon_k^\ast = \epsilon_{k-1}^\ast(\mu_k)$
	\If{$\varepsilon_k^\ast \leq \epsilon_{\text{tol, min}}$}
		\State $N_{\max} = k - 1$
		\State exit
	\EndIf
	\State $S_k = S_k \cup \mu_k$
	\State $W_k = \text{span} \{W_k, \{ u(\mu_k)\}\}$
\EndFor
\end{algorithmic}
\end{algorithm}

Recall that the set of basis functions $\{ u(\mu_1), u(\mu_2), \ldots, u(\mu_N)\}$ are usually pointing in similar directions, so in order to make (\ref{eq:galerkin}) well-conditioned, the basis functions need to be orthogonalized. 
This process can be incorporated into Algorithm~\ref{alg:ga} for efficiency. 
Note that this algorithm does not give a unique set $W_N$, but the error given by $\varepsilon_k^\ast$ is monotone decreasing~\cite{binev}. 

\subsection{Gradient Algorithm}\label{sec:grad}
We can consider the choice of an appropriate sample set $S_N$ by describing the reduced basis choice as an minimization problem of choosing the best generating set $\underline{\mu}= \lbrace \mu_1,\mu_2,\ldots , \mu_N\rbrace$, for $\mu_i\in \mathcal{D}$ for $1\leq i\leq N$ and a fixed $N\in\mathbb{N}$. Based on this idea of minimizing a function we use a gradient algorithm to find the desired sample $S_N$ by solving the minimization problem 
\begin{equation}
\mu_k = \underset{\mu \in D}{\text{arg min}} \sum_{\nmu \in \Upsilon} \frac{\Vert u^{fe}(\nmu)-u_{N+1}(\nmu)\Vert_X}{\Vert u^{fe}(\nmu)\Vert_X}
\label{eq:min1}
\end{equation}
in each iteration. In this minimization problem, $u_{N+1}(\nmu)$ describes the reduced basis approximation for one additional parameter in $S_{N+1}= S_N\cup \mu$. $\Upsilon \subset \mathcal{D}$ denotes a mesh over the parameter space. 
 Thus we minimize a function of $\mu$. We choose the first sample $\mu_1$ randomly, analogously to the greedy algorithm. Next, the sample set $S_N$ is constructed one additional sample at the time by solving the minimization problem above.

By finding $\mu_k$ over $\mathcal{D}$, instead of over a fine mesh $\Xi\subset \mathcal{D}$ as  in the greedy algorithm, we can improve the sample set $S_N$ and therefore reduce its dimension, $N$. This depends on the chosen \textit{a proiri} error bound
\[\Vert u^{fe}(\mu)-u_N(\mu)\Vert_X\leq\epsilon_{tol}\] 
and the maximum dimension $N_{max}$. Additionally, a very fine mesh  $\Xi\subset \mathcal{D}$ is needed to ensure an appropriate sample set $S_N$ using the greedy algorithm. This yields a large amount of calculations for each step. Using iterative methods to solve the minimization problem, we can reduced the required calculations significantly. However, in order to ensure convergence of the algorithms, more requirements must be meet. Below, we discuss these requirements further.

For all iterative methods, in order to determine a local minimum a starting point has to be chosen. To avoid clustering of the $\mu_i\in \mathcal{D}$ around a certain parameter $\xi$ due to local minima, alternating the starting point in each iteration is advised. First, a coarse mesh $\Lambda$ of the parameter space $\mathcal{D}$ is created. 
 From this mesh a starting point $\mu_k^0$ is chosen by  
\[\mu_k^0 = \underset{\mu \in \Lambda}{\text{arg min}} \sum_{\nmu \in \Upsilon} \frac{\Vert u^{fe}(\nmu)-u_{N+1}(\nmu)\Vert_X}{\Vert u^{fe}(\nmu)\Vert_X}.\]
This way we reduce the likelihood of clustering $\mu_i\in \mathcal{D}$. In addition, we address the problem regarding local instead of global minima. For nonlinear functions, numerical algorithms find mostly local minima~\cite{opti1, opti2}. By choosing the starting point from a coarse mesh over the parameter space we can improve the determined local minimum. This approach is similar to the hill climbing algorithm~\cite{rus10}. However, there is no guarantee that the determined local minimum for the starting point $\mu_k^0$ is indeed a global minimum over $\mathcal{D}$. 
To find the desired local minimum we use a gradient method, as it most accurately provides the desired minima without too many requirements on the function. Alternatively the Nealder-Mead algorithm (also known as the simplex algorithm) can be used if not enough of the requirements for a gradient method are met~\cite{alt02}. 

Both the gradient and Nealder-Mead algorithms are non-restricted methods. Adding a penalty term $\phi(\mu)$ to the minimization functional
 converts the restricted minimum problem to a non-restricted minimum problem. To ensure convergence using a penalty term, the parameter space, $\mathcal{D}$, needs to be closed and compact. Moreover, the minimization function needs to be in $C^1(\Omega)$ which means the derivative exists. The standard gradient method does not use an efficient step size for each iteration, as it only defines the direction of steepest decent. We use the Armijo-Powell rule to find an efficient step size and thus ensure convergence of the algorithm. 

Next, we introduce two stopping criteria. The first one is an \textit{a priori} maximum number of elements, $N_{\max}\in \mathbb{N}$, in the sample set. The second one is an \textit{a proiri} error bound, as described above. These stopping criteria are the same as for the greedy algorithm. 

Overall, we used nonlinear optimization methods to improve the greedy approach to finding a suitable sample set for the reduced basis approximation. While this may improve calculation speed and accuracy, it can be difficult to show that the minimization problem meets all requirements. This depends on the PDE at hand.
A summary of this method is given by Algorithm~\ref{alg:gra}.
 
  \begin{algorithm}
\caption{Gradient Sampling Algorithm}\label{alg:gra}
\begin{algorithmic}
\State Initialize Setup
\State $\mu_1 = rand(\Xi\subset \mathcal{D})$
\For{$k = 2:N_\text{max}$}  
\State Find starting point for minimization 
\State $\mu_k^0 = \underset{\mu \in \Lambda}{\text{arg min}} \sum_{\nmu \in \Upsilon} \frac{\Vert u^{fe}(\nmu)-u_{N+1}(\nmu)\Vert_X}{\Vert u^{fe}(\nmu)\Vert_X}$
\State Solve minimization problem
\State $\mu_k = \underset{\mu \in D}{\text{arg min}} \sum_{\nmu \in \Upsilon} \frac{\Vert u^{fe}(\nmu)-u_{N+1}(\nmu)\Vert_X}{\Vert u^{fe}(\nmu)\Vert_X}$
\State $\varepsilon_k= \max $
\If{$\varepsilon_k^\ast \leq \epsilon_{\text{tol, min}}$}
\State $N_{\max} = k - 1$
\State exit
\EndIf
\State $S_k = S_k \cup \mu_k$
\State $W_k = \text{span} \{W_k, \{ u(\mu_k)\}\}$
\EndFor
\end{algorithmic}
\end{algorithm}

Further improvements on the convergent rate can be achieved using more complex methods such as conditional gradient and SQP methods. 
Similar to the comparison of the gradient algorithm to the Nealder-Mead algorithm, these improvements come at a cost. In our case the cost is the additional requirements on the minimization function. For more information about nonlinear optimization in general, and conditional gradient as well as SQP methods, refer to~\cite{alt02}.

\subsubsection{A few remarks on implementation.}
\label{sec:gradimp}
 In this work, we implemented an additional change to the gradient algorithm. To be able to compare the gradient algorithm with the other methods (specifically the greedy algorithm), a larger basis size is taken into consideration even if the stopping criteria was already met. We noted that in some cases, depending on the choice of $\Lambda$, the gradient method performed matrix operations on ill-conditioned matrices and began to become instable. A sensitivity analysis showed that this was due to clustering of the $\mu_k^0$ for different $k$s, which was also observed using the greedy algorithm. 

Two solutions to this problem are possible. To avoid clustering for large $N$ the start of the minimization iteration has to be different from the last iteration. That is,
\[\mu_{k-1}^0\neq \mu_k^0.\]
To ensure higher stability, the conditioning of the matrix can be improved by using Gram-Schmidt orthogonalization within the gradient sampling algorithm. This increases the computation time for each iteration, but the matrices involved in the computation are no longer ill-conditioned. In this work a combination of both solutions was used, but further tests showed that the orthogonalization would have been sufficient to solve the problem.

\subsection{Metropolis Algorithm} 
In the previous section, we described the construction of a reduced basis
 as the minimization problem of choosing the best basis generating set $\underline{\mu}=\{ \mu_1, \mu_2, ..., \mu_N \}$, for $\mu_{i}\in \mathcal{D}, 1\leq i\leq N$ and a fixed $N\in\mathbb{N}$. Note that $\underline{\mu} = S_N$ for fixed $N$. In this section, we will consider an alternative to the minimization problem~(\ref{eq:min1}), given by
\begin{equation}\label{minprob}
\min_{\underline{\mu}\in \mathcal{D}^N} \sum_{\mu\in \Upsilon}\frac{||u^{fe}(\mu)-u_N(\mu;\underline{\mu})||_{X}}{||u^{fe}(\mu)||_{X}}.
\end{equation}
The goal in this section will be to find the Bayes' estimate of $\{ \tilde{\mu}_1, \tilde{\mu}_2, ..., \tilde{\mu}_N \}$, where $\tilde{\mu}_{i}$ denotes the random variable of $\mu_i$.
Note that, contrary to the methods already described, the minimization problem in (\ref{minprob}) solves for the entire basis generating set at once, rather than for one element at a time. Thus, the algorithm must be started over for each new $N$ to find a new basis set $S_N$. 
In order to work on finding the Bayes' estimate, we first need to find the posterior probability density function of $\{ \tilde{\mu}_1, \tilde{\mu}_2, ..., \tilde{\mu}_N \}$.

\subsubsection{The Posterior Probability Density Function }

We will use Bayes' theorem to obtain the  posterior probability density function of $\tilde{\mu}_1, \tilde{\mu}_2, ..., \tilde{\mu}_N $ given the error $\delta=\left\{u(\mu) - u^{fe}(\mu)\right\}_{\mu\in \Upsilon}$. If we denote $\tilde{\delta}$ as the random variable  of $\delta$, then the posterior probability density becomes
\begin{equation}\label{posterior}
P_{\tilde{\mu}_1, \tilde{\mu}_2, ..., \tilde{\mu}_N}(\underline{\mu} | \delta)\propto P_{\tilde{\delta}}( \delta|\underline{\mu})P_{\tilde{\mu}_1, \tilde{\mu}_2, ..., \tilde{\mu}_N}(\underline{\mu}),
\end{equation} 
where $P_{\tilde{\delta}}( \delta|\underline{\mu})$ represents the density of the noise and $P_{\tilde{\mu}_1, \tilde{\mu}_2, ..., \tilde{\mu}_N}(\underline{\mu})$ the prior density of $\tilde{\mu}_1, \tilde{\mu}_2, ..., \tilde{\mu}_N$. Considering the minimization problem (\ref{minprob}), the natural choice for the density of $\tilde{\delta}$ given $\underline{\mu}$ is
\begin{equation}
P_{\tilde{\delta}}( \delta|\underline{\mu})\propto \exp\left(- \sum_{\mu\in \Upsilon}\frac{||u^{fe}(\mu)-u_N(\mu;\underline{\mu})||_{X}}{||u^{fe}(\mu)||_{X}} \right).
\end{equation} 
There are several reasonable choices for the prior density. In general we \textbf{can} include any prior knowledge of the properties of $\tilde{\mu}_1, \tilde{\mu}_2, ..., \tilde{\mu}_N$.
However, to allow for a fair comparison between the Greedy algorithm and the MCMC method, we will restrict ourselves to the most basic prior knowledge. Particularly, we assume that the parameters are in the parameter space and that the parameters are sorted in ascending order, i.e. $\mu_1, \mu_2, ..., \mu_N \in \mathcal{D}$ and $\mu_1 < \mu_2 < ... <\mu_N$. This prior density can be expressed as
\begin{equation}
P_{\tilde{\mu}_1, \tilde{\mu}_2, ..., \tilde{\mu}_N}(\underline{\mu})\propto \chi_{\mathcal{E}}(\mu_1, \mu_2, ..., \mu_N)=\left\{
\begin{array}{l}
1 \text{ if } \mu_1, \mu_2, ..., \mu_N\in \mathcal{E}\\
0 \text{ else}
\end{array}
\right. \label{equKE}
, 
\end{equation} 
where $\mathcal{E}=\{\mu_1, \mu_2, ..., \mu_N\in \mathcal{D}| \mu_1 < \mu_2 < ... <\mu_N\}$, and $\chi_{\mathcal{E}}$ is the indicator function over the set $\mathcal{E}$.

By definition, the solution of the statistical inverse problem is the posterior density for the $\tilde{\mu}_1, \tilde{\mu}_2, ..., \tilde{\mu}_N$, but we are interested in finding the Bayes' estimate $E(\tilde{\mu}_1, \tilde{\mu}_2, ..., \tilde{\mu}_N| \delta)$. Unfortunately, the posterior density (\ref{posterior}) does not have a closed form, which makes it impossible to find the Bayes' estimate in a direct way.
Therefore, we approximate the Bayes' estimate via simulation.

\subsubsection{The Markov Chain Monte Carlo Method}

The idea of the MCMC method is to generate a large random sample from the posterior density $P_{\tilde{\mu}_1, \tilde{\mu}_2, ..., \tilde{\mu}_N}(\underline{\mu} | \delta)$, with Y samples, and then approximate the Bayes' estimate by the sample mean,
\begin{eqnarray} \label{estEqu}
E(\underline{\mu} | \delta)=\int_{\mathbb{R}^N} \underline{\mu} P_{\tilde{\mu}_1, \tilde{\mu}_2, ..., \tilde{\mu}_N}(\underline{\mu}| \delta) d\underline{\mu} \approx  \frac{1}{Y}\sum_{i=1}^{Y} \underline{\mu}_i,
\end{eqnarray}
where $ \underline{\mu}_i$ represents the $i^{\text{th}}$ random sample of the posterior density (\ref{posterior}). Typical algorithms to generate such large random samples are the Gibbs sampler or the Metropolis Hastings algorithm~\cite{Chib1995}. In this manuscript we will use a Metropolis Hastings algorithm and mention some adaptive versions of it. 

Algorithm \ref{PilotA} describes the Metropolis Hastings algorithm. Let $q(\underline{\mu}_{i-1},\underline{\mu}_{*})$ be a candidate-generating density, that is, a density which is used to generates a new candidate random sample $\underline{\mu}_{*}$ from a current random sample $\underline{\mu}_{i-1}$. We assume that $q_{C}(\underline{\mu}_{i-1},\cdot)$ is a density function with mean $\underline{\mu}_{i-1}$ and covariance matrix $C$.  We choose $q_{C}(\underline{\mu}_{i-1},\cdot)$ to represent a Gaussian density. Assume that $\underline{\mu}_{*}$ is a generated candidate then it will be accepted as a true random sample $\underline{\mu}_{i}$ with probability 
\begin{equation} \label{alpha}
\alpha(\underline{\mu}_{i-1},\underline{\mu}_{*})=\left\{
\begin{array}{c l}
    \min \left[\frac{P(\underline{\mu}_{*}|\delta)q(\underline{\mu}_{*},\underline{\mu}_{i-1})}{P(\underline{\mu}_{i-1}|\delta)q(\underline{\mu}_{i-1},\underline{\mu}_{*})}, 1 \right], & \text{if }  P(\underline{\mu}_{i-1}|\delta)q(\underline{\mu}_{i-1},\underline{\mu}_{*})>0 \\
    1, &  \text{otherwise},
\end{array}\right. 
\end{equation}
where $P(\cdot|\delta)$ represents the posterior density. We refer to  $\alpha(\underline{\mu}_{i},\underline{\mu}_{*})$ as the acceptance ratio. If $q_{C}(\underline{\mu}_{i-1},\underline{\mu}_{*})$ is a symmetric density function, (\ref{alpha}) can be simplified by $q(\underline{\mu}_{*},\underline{\mu}_{i-1})=q(\underline{\mu}_{i-1},\underline{\mu}_{*})$.

Note that Metropolis Hastings algorithm only generates true random samples from the posterior density after a Burn-in time $B$. Hence, the first $B$ random samples must be ignored when using~(\ref{estEqu}) to obtain the Bayes' estimate. 

\begin{algorithm}
\caption{The Metropolis Hastings Algorithm}
\label{PilotA}
\begin{algorithmic}
    \State Take any $\underline{\mu}_{0}$ from $\mathcal{E};$
    \For{i = 1 to B+Y}
    \State  Generate $\underline{\mu}_{*}$ from $q_{C}(\underline{\mu}_{i-1},\cdot)$ and $u$ from $\mathcal{U}(0,1)$;
     \If{$u \leq \alpha(\underline{\mu}_{i},\underline{\mu}_{*})$}
     	\State $\underline{\mu}_{i}=\underline{\mu}_{*};$
	\Else
	\State  $\underline{\mu}_{i}=\underline{\mu}_{i-1};$
	\EndIf
     \EndFor\\

    \State   \textbf{Return} $\left\{ \underline{\mu}_1, \underline{\mu}_2, ..., \underline{\mu}_{B}, \underline{\mu}_{B+1}, ..., \underline{\mu}_{B+Y}\right\}$
  \end{algorithmic}
\end{algorithm}

It is known that the proper choice of the proposed density $q_C(\cdot, \cdot)$ for the Metropolis algorithm is vital to obtain a reasonable result by simulation in a suitable amount of time. This choice is generally very difficult since the target density is generally unknown \cite{Haario2001, Gelman1996, Gilks1998}. One possible way of smoothing out this problem is by using an adaptive Metropolis algorithm which iteratively updates the proposal density (or its covariance matrix) in an appropriate way. The downside of this kind of adaptive algorithm is that the chain usually becomes non-Markovian, which would require establishment of the correct ergodic properties. In this manuscript we used the algorithm proposed in \cite{thilo} to find a proper covariance matrix for the candidate-generating density $q_{C}(\cdot,\cdot)$. However, there are many other adaptive Metropolis algorithms, such as \cite{Gilks1998, Haario2001} to name a few.

We saw that the Metropolis-Hastings algorithm samples properly from the posterior distribution only after a burn-in time $B$. Therefore, we are interested in knowing how long this burn-in phase should be. Unfortunately, there is no theory giving a good estimate of the burn-in time prior to running the Metropolis-Hastings algorithm. That leaves us to first run the Metropolis-Hastings algorithm and then check whether or not it converged to its invariant distribution. There are several methods to check whether or not the chain converged, however all methods have a positive probability of producing the wrong conclusion. Therefore, it is recommended to use more than one diagnostic method. The most common diagnostic methods in the literature are found in the work of Gelman and Rubin \cite{Gelman1992,Brooks1998}, Geweke \cite{Adlouni2006} and Raftery and Lewis \cite{Raftery1995}.

 \section[H]{Hyperspectral Diffuse Optical Tomography Model}
Optical tomography (OT) is an imaging technique that images an object through sectioning by use of an optical wave, that is, a wave from a light source. 
Typically, this source is laser light in the visible (about 400 to 700nm) or near infrared range (about 700 to 1600 nm).
Tissue is a highly scattering medium and so, as the collimated laser beam passes through the tissue some of the light is absorbed by chromophores (such as hemoglobin, lipid and water), but most is scattered. 
Detectors placed on the boundary of the tissue collect the scattered beams, and from this data a 2-D image (slice) of the tissue is reconstructed in the form of a spatial map of the tissue's absorption and scattering coefficients~\cite{jiang}. 
Both of these coefficients are dependent on the wavelength, $\lambda$, of the light source~\cite{durduran}. 
Since cells in tumors have higher absorption coefficients than normal cells due to an increased water or ionic concentration, and they also scatter photons differently, the absorption and scattering coefficients of the cells being imaged are the most important parameters to be determined in most medical applications~\cite{jiang}. 

The relationship between the density of the light source at each location, $\mathbf{x}$, in the tissue, and the absorption and scattering coefficients is typically described using 
 the diffusion approximation to the radiative transport equation, which results in a modality known as diffuse optical tomography (DOT). 
In the frequency domain, the forward problem is expressed using the photon diffusion model 
\begin{align}
-\nabla \cdot (D(\mathbf{x})\nabla u(\mathbf{x})) + (\mu_a(\mathbf{x}) + ik)u(\mathbf{x})& = h(\mathbf{x})& \text{in } \Omega
\label{eq:DOTforward}\\
-D \frac{\partial u}{\partial \mathbf{n}} & = f & \text{on } & \partial \Omega, 
\label{eq:DOTforward3}
\end{align}
an elliptic PDE, where $u$ is the photon density, $h$ is an interior forcing function, $k$ is the wave number of the modulating frequency of the laser, and $D$ is the diffusion coefficient, expressed as $D = 1/3 (\mu_a + \mu_s')$ where $\mu_a$ is the absorption coefficient and $\mu_s'$ is the reduced scattering coefficient~\cite{arridgereview, arridgereview_new, jiang}. 
The diffusion model is a first-order approximation to the radiative transport equation, assuming $\mu_s' \gg \mu_a$ and the detector and source are not too close together~\cite{arridgereview_new, jiang}.

Hyperspectral imaging (HSI) is an imaging modality involving the use of hundreds of optical wavelengths that has been generating a lot of research in recent years, specifically regarding its application to remote sensing, and geospatial imaging~\cite{chen_sparse,shaw,sun}. 
Currently, researchers are looking at HSI's applications in the medical field as a tool for non-invasive disease diagnosis and surgical guidance~\cite{fei_hyperspect,kim,lufei,poole14}. 
HSI is a hybrid modality, combining imaging with spectroscopy by collecting spectral information at each pixel of an array while scanning through different wavelengths to generate a three-dimensional hypercube of spatial and spectral data~\cite{lufei}. 

It has recently been proposed to apply hyperspectral imaging to DOT to increase the spectral resolution of the image while maintaining the high spatial resolution given by the optical imaging modality. 
Application of HSI to DOT results in what is known as hyperspectral DOT (hyDOT).
HyDOT has an added dimension as the governing PDE given in (\ref{eq:DOTforward}) is now parameterized by the wavelength, $\lambda$. 
\begin{equation}
-\nabla \cdot (D(\mathbf{x},\lambda)\nabla u) + (\mu_a(\mathbf{x},\lambda) + i\kappa)u= h(\mathbf{x})
\label{eq:hyDOTforward}
\end{equation}
Note that in the general discussion of the RBM in the previous section the parameter of a parameterized PDE was denoted as $\mu$ but, in our specific application, as described in Section~\ref{compare}, the parameter is the wavelength, $\lambda$. 
The main goal in hyDOT is to solve the image reconstruction problem, which is also known as the inverse problem. 
That is, given the data collected by the sensors on the boundary, reconstruct a spatial map of the absorption ($\mu_a$) and scattering coefficients ($\mu_s'$ where $D = \mu_s' + \mu_a$) for the tissue. 
This problem is ill-posed and is usually solved numerically using an iterative method where the forward problem is solved hundreds of times. 
This paper focuses on reducing the computational burden of the forward problem with the ultimate goal of solving the inverse problem more efficiently. 

\section{Comparison of the Algorithms Using Simulation}
\label{compare}
 \subsection{Sample Problem}
\begin{figure}[h]
\centering
\begin{tikzpicture}
 \draw[left color=darkblue, right color=darkblue!40!white,shading angle=135, drop shadow] (0,0) circle (2.5);
 \draw[left color=myred, right color=myred!50!white,shading angle=135] (-1,-1) circle (0.8);
 \draw[thick, ->, shadowed] (0,-3) -- (0,3);
 \draw[thick, ->, shadowed] (-3,0) -- (3,0);
 \node at (-1,-1) {\Large{$\Omega_1$}};
 \node at (1,1) {\Large{$\Omega_0$}};
\end{tikzpicture}
\caption[Geometry of the simple example for application of RBM to the forward problem of hyDOT]{The layout of geometry of the simple problem. We have a 2-D circular sample of tissue of radius 25 centimeters, centered at the origin, with a circular tumor of radius 5 centimeters located at the point (-15,-10). The domain is split into the healthy tissue, $\Omega_0$ and the cancerous cells, $\Omega_1$. }
\label{fig:simpgeo}
\end{figure}
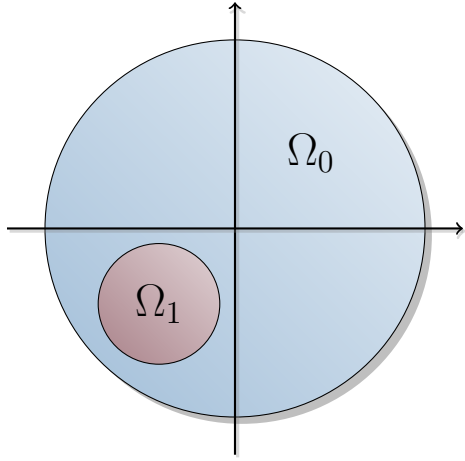

 To illustrate the application of this method, we consider a simplified version of the forward problem for hyDOT where $\kappa$, the wave number for the modulation of the laser light 
and $h$, the interior forcing function, are considered to be zero.
We consider only Neumann boundary conditions with one Gaussian source given by 
\[f(x,y) = 15e^{\frac{-\left((x - x_1)^2 +(y-y_1)^2\right)}{10}} \]
where $(x_1, y_1) \approx (-24.5196, -4.8773)$, one of the finite element mesh points on the boundary.
The strong form of the governing PDE is given by
\begin{eqnarray}
- \nabla \cdot (D(\mathbf{x},\lambda)\nabla u) + \mu_a(\mathbf{x},\lambda) u &=& 0 \qquad \text{in } \Omega \label{eq:simpstrong}\\
D(\mathbf{x},\lambda) \frac{\partial u}{\partial \nu} &=& f \qquad \text{on } \partial \Omega \label{eq:simpstrongb1}
\end{eqnarray}
The wavelength, $\lambda$, is the parameter on which the PDE is dependent.  
We will consider the parameter space as $\mathcal{D} = [600, 1000] \subseteq \mathbb{R}$, with units nanometers, which is a typical range used in practice.
We consider a simple geometry for $\Omega$ in two dimensions given by a circle of radius 25 centimeters centered at the origin on a Cartesian grid, with a circular tumor of radius $5$ centimeters located at the point (-15, -10)(see Figure~\ref{fig:simpgeo}).
The domain, $\Omega$ is split into two regions, namely the healthy tissue, $\Omega_0$ and the cancerous cells, $\Omega_1$.
The location of the tumor was chosen to be relatively close to the source to yield results that were easier to visualize, since light does not penetrate very far into tissue. 
In the forward problem, the source, $f$, and the geometry and location of the area of interest, in our case a collection of cancerous cells, is known. 
Only the measurements of the scattered photons on the boundary are unknown and are given by the solution, $u$. 

We wish to use the RBM to find an approximate solution $u_N(x, \lambda)$ to (\ref{eq:simpstrong}) - (\ref{eq:simpstrongb1}). 
The weak formulation of the PDE is given by
\begin{equation} 
\int_{\Omega} \left(D\nabla u \cdot \nabla v + \mu_auv\right)dx =  \int_{\partial \Omega} f \gamma_D( v) ds \qquad \forall v \in V 
\label{eq:weaksimp}
\end{equation}
where $V$ is the space of smooth test functions, and $\gamma_D: H^1(\Omega) \rightarrow H^{1/2}(\partial \Omega)$, the Dirichlet trace. 
From the weak formulation, we define the following bilinear and linear forms, 
\begin{equation}
a(u,v;\lambda) = \int_\Omega (D(\lambda)\nabla u \cdot \nabla v + \mu_a(\lambda)uv) dx  \qquad F(v) = \int_{\partial \Omega} f\gamma_D(v) ds 
\label{eq:simpbilinear}
\end{equation} 
where $a(u,v): H^1(\Omega) \times H^1(\Omega) \to \mathbb{R}$ and $F(v): H^1(\Omega) \to \mathbb{R}$, with the inner product induced by the bilinear form. 
Note that here we are in compliance (that is, the output functional $\ell = f$). 

\begin{figure}[h]
\begin{center}
	\includegraphics[width=0.7\textwidth]{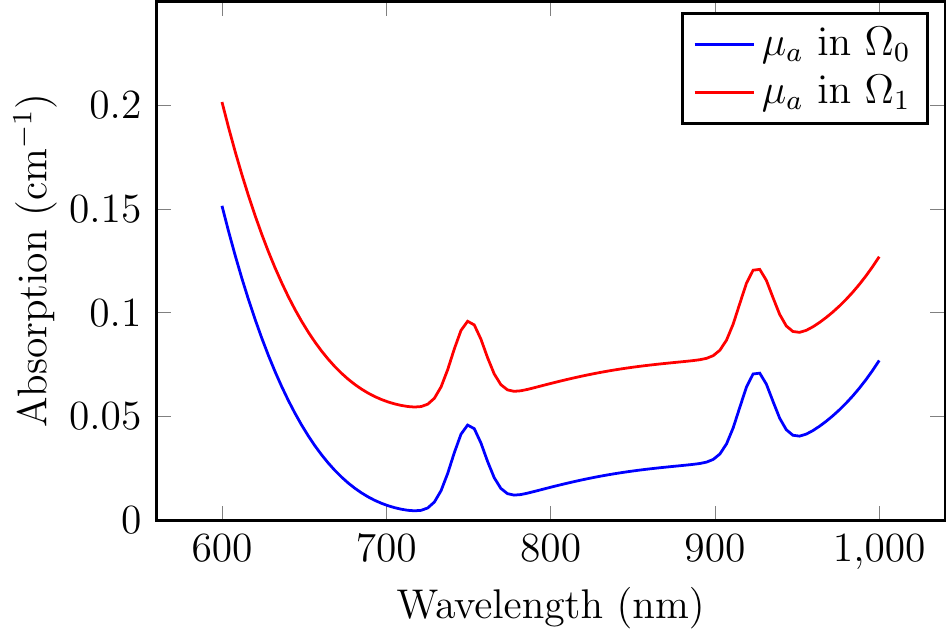}
\end{center}
\caption[Absorption coefficient as a function of $\lambda$, simple example]{Graphs of the absorption coefficient, $\mu_a$, (in $\text{cm}^{-1}$) as a function of wavelength, $\lambda$, (in nm) in the healthy tissue, $\Omega_0$ (blue) and in the cancerous tissue, $\Omega_1$(red).}
\label{fig:mua}
\end{figure}

We created functions for $\mu_a(\lambda)$ such that given $\lambda$ the value of $\mu_a$ would be one constant value anywhere in $\Omega_0$ and a different constant value anywhere in $\Omega_1$. 
Since 
\[ D = \frac{1}{3(\mu_a + \mu_s')} \]
the value for the diffusion coefficient would follow a similar pattern as it is a function of $\mu_a$. 
The function values were based on the graph showing the absorption coefficient as a function of wavelength given by Saibaba et al. in~\cite{saibaba} and experimental values given by Yodh and Chance~\cite{yodh} and Jiang~\cite{jiang}. 
Given these experimental values, the function for $\mu_a$ in $\Omega_0$ was taken to be a quartic check function (found using interpolation through select data points) with Gaussian spikes at 725 and 950 nanometers. 
Since cancerous cells generally have higher absorption coefficients than healthy tissue,
we followed the example of \cite{saibaba} and made the profile of the absorption coefficient in $\Omega_1$ to be a positive perturbation of its profile in $\Omega_0$. 
Both profiles are given in Figure~\ref{fig:mua}.
The reduced scattering coefficient $\mu_s'$ was chosen to be 17 $\text{cm}^{-1}$ at all points in the spatial domain. 
This value was chosen because Yodh and Chance~\cite{yodh} reported that at 820 nm, they experimentally obtained $\mu_s'$ values of 16.5-18.5 $\text{cm}^{-1}$ in human brain tissue and values of 12.7- 17.3 $\text{cm}^{-1}$ in human breast tissue. 

To choose the set of fixed, nested parameter samples, $S_n = \{ \lambda_j,\}_{j=1}^n$ where $n \ll N$, we compare Algorithms~\ref{PilotA},\ref{alg:gra} to the greedy algorithm given in Algorithm~\ref{alg:ga} where $\epsilon_{k}^\ast$ is the dual norm of the residual as given in~(\ref{eq:dualnorm}), with the norm on $X$ the standard $H^1$ norm. 
The test space, $\Xi$, was generated using a fine linear mesh of $M=400$ equally spaced points between 600 and 1000 nanometers.
The tolerance was set at 1e-5 for the greedy algorithm and 1e-7 for the gradient algorithm (to force the algorithm to choose a basis of a size $N >4$). 
We let $N_{\max}$ vary for bases of sizes in the set $\{5,6,7,8,9,10,15,20\}$.

We constructed and stored solutions for each value of the parameter in $S_N$ using the finite element discretization, 
\[\tilde{u}_j = \tilde{u}(\lambda_j) = \sum_{i=1}^\mathcal{N} b_i(\lambda_j)\phi_i. \] 
The finite element solution was found using the PDE Toolbox in Matlab over a mesh of 2097 elements~\cite{matlab}. 
Then we defined an approximation space $W_N = \text{span}\{\hat{u}_j \}_{j=1}^N$. 
Note that in both cases, $W_N$ is constructed as a Lagrange approximation subspace since the derivatives of $u$ with respect to $\lambda$ are not known.  
The simulations were run on a MacBook Pro, OS X Version 10.9.5, with a 2.5 GHz Intel Core i5 processor, and 4 GB 1600 MHz memory. 
We note, however, that the gradient algorithm was later updated to improve stability, as described in Section~\ref{sec:gradimp}, and the updated algorithm was run on a Mac, OS X Version 10.12.6, with a 3.2 GHz Intel Core i5 processor, and 32 GB 1867 MHz memory,
as the previous machine was no longer available. 
The relative errors were comparable to previous results, though the running time was slightly faster, as we would expect with a faster processor. That said, the run times on the new machine were still within the same order of magnitude, so we maintain that the results from the different computers are comparable. 

Since the efficiency of the RBM relies on the affine decomposition of the bilinear form, we must demonstrate that $a(u,v; \lambda)$ given in (\ref{eq:simpbilinear}) has affine parameter dependence. 
We first note that, although $\lambda$ does not have spatial dependence, the diffusion coefficient $D$ does. 
Thus, we must first decompose the bilinear form geometrically, considering the domain $\Omega_0 = \{(x_1,x_2) | x_1^2 + x_2^2 \leq 625\} \backslash \Omega_1$ of the healthy tissue, and the domain $\Omega_1 =\{ (x_1,x_2) |  (x_1+15)^2 + (x_2+10)^2 \leq 25\}$ of the cancerous tissue. 
Since the diffusion and absorption are homogeneous within each of these domains by construction, we can consider the functions, $D_0$, $\mu_a^0$ and $D_1$, $\mu_a^1$ on $\Omega_0$ and $\Omega_1$, respectively, that are functions of $\lambda$ only (that is, they are spatially independent). 
Therefore, we can decompose the bilinear form as 
\begin{eqnarray*}
a(u,v;\lambda) &=& \sum_{q=1}^Q \Theta^q(\lambda)a^q(u,v)\\
&=& D_0(\lambda) \int_{\Omega_0} \nabla u \cdot \nabla v dx + \mu_a^0(\lambda) \int_{\Omega_0} uv dx\\
&+& D_1(\lambda) \int_{\Omega_1}\nabla u \cdot \nabla v dx + \mu_a^1(\lambda) \int_{\Omega_1}uvdx.
\end{eqnarray*}
In this case, the linear form $F$ has no explicit parameter dependence and so it does not need to be decomposed. 
The reduced basis approximation to the problem is then computed by solving the problem
\[ \hat{A}_\lambda \hat{C} = \hat{F}\] 
where $\hat{C}$ is the vector of unknown coefficients, 
 and $\hat{A}_\lambda = C^TA_\lambda C$ where 
\begin{eqnarray}
A_\lambda &=& D_0(\lambda)A_{00} + \mu_a^0(\lambda)A_{01} + D_1(\lambda)A_{10} + \mu_a^1(\lambda)A_{11} \\
&(A_{00})_{i,j}& = \int_{\Omega_0} \nabla u_i(x) \cdot \nabla u_j(x) dx \qquad 1 \leq i,j, \leq N\\ 
&(A_{01})_{i,j}& = \int_{\Omega_0}  u_i(x) \cdot u_j(x) dx \qquad 1 \leq i,j, \leq N\\ 
&(A_{10})_{i,j}&=\int_{\Omega_1} \nabla u_i(x) \cdot \nabla u_j(x) dx \qquad 1 \leq i,j, \leq N\\
&(A_{11})_{i,j}& = \int_{\Omega_0}  u_i(x) \cdot u_j(x) dx \qquad 1 \leq i,j, \leq N
\end{eqnarray}
and $C$ contains the coefficients for the finite element approximation, $u_j$, for $j = 1, \ldots, N$. 
Due to the decomposition of $A$, 
the matrices $A_{00}, A_{01}, A_{10}, A_{11}$
can be precomputed and stored in the offline stage. 
We note that orthogonalization of the basis functions using the Gram-Schmidt method with respect to the inner product induced by the bilinear form was necessary for the conditioning of the matrix $A$.
Without orthogonalization the condition number of the matrix was as high as order $10^{20}$ whereas with orthogonalization the condition number was generally around 1, and no higher than 10 at select wavelengths. 

For problems in which the bilinear form $a(u,v;\mu)$ has an affine decomposition of the form, 
\[ a(u,v;\mu) = a_0(u,v) + \mu a_1(u,v) \] 
where $a_0(u,v)$ is symmetric, coercive and continuous, and the parameter space is given by $D = [0, \mu_{\text{max}}] \subset \mathbb{R}$, 
it has been shown that if the sample points are chosen according to a logarithmic point distribution, 
then the reduced basis approximation will converge exponentially to the truth solution for all $N$ greater than some value $N_{\text{crit}}$~\cite{maday}. 
Further, it has been numerically demonstrated that using a logarithmic versus a uniform or Chebyshev distribution in the choosing of points results in a much smaller maximum relative error~\cite{veroy_thesis}. 
Thus, we decided to compare all three basis-selection methods to a logarithmic distribution of points, where each of the $\lambda_n$ satisfied
\[ \ln(\bar{\sigma}\lambda_{n+1}) = \frac{n-1}{N-1} \ln(\bar{\sigma}\lambda_{\max} + 1) \]
where $\lambda_{\max} = 1000$ nm, and $\bar{\sigma}$ is an upper bound for the maximum eigenvalue of $a_1(u,v)$ relative to $a_0(u,v)$, which we set as $5.5\times 10^4$. 

\subsection{Results}
\begin{figure}[h] 
\begin{center}
\includegraphics[width=\textwidth]{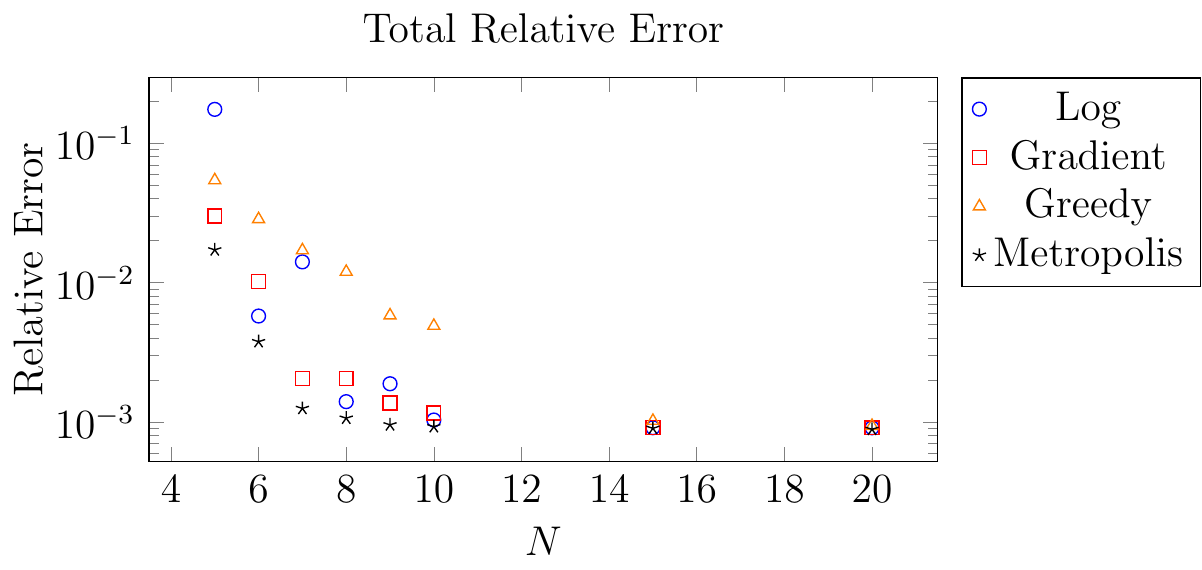}
\end{center}
\caption[Total relative error for all three methods]{The total relative error over 100 wavelengths for the solutions calculated using the bases generated by the greedy, Metropolis, gradient, and log-scale algorithms, respectively, as a function of the size of the basis, $N$. An average of ten trials was used for the Metropolis and greedy methods.}
\label{fig:relerror}
\end{figure}

\begin{figure}[h] 
\begin{center}
\includegraphics[width=\textwidth]{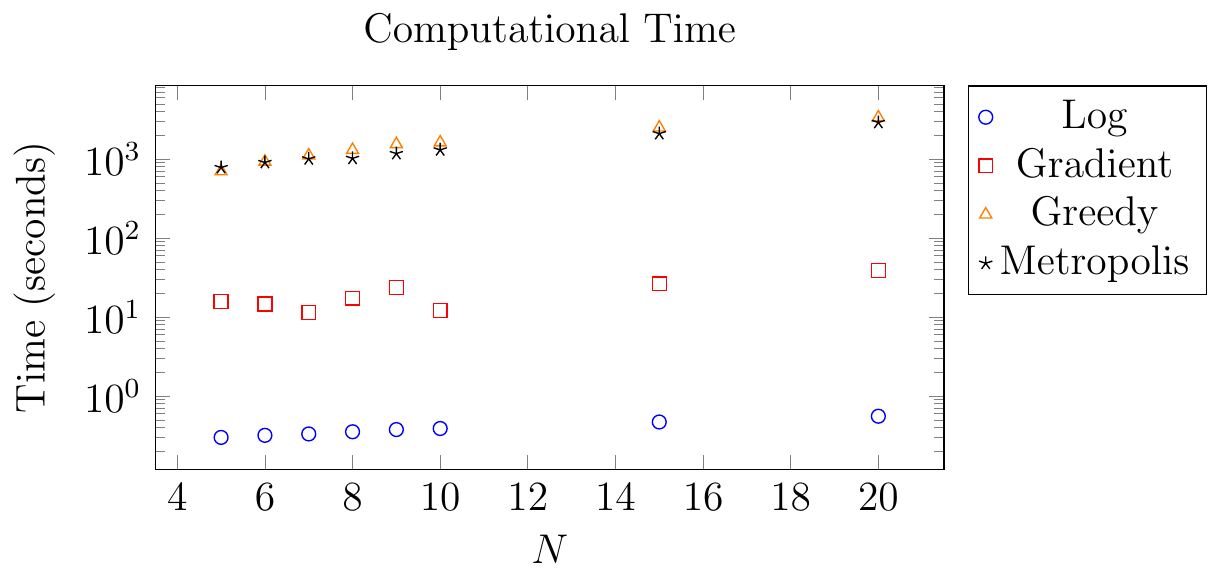}
\end{center}
\caption[Timing for all three methods]{The computational time in seconds to find the wavelengths that will generate the basis for the greedy, Metropolis, and gradient, respectively, as a function of the size of the basis, $N$. An average of ten trials was used for the Metropolis and greedy methods. Note that for the log-scale method there is minimal to no time involved in finding the basis.}
\label{fig:timing}
\end{figure}

To test the efficiency of the RBM on this simple example, given a basis generated by one of the three algorithms described above (in addition to using the logarithmic spacing), the reduced basis approximation of the solution was computed for 100 linearly spaced values of $\lambda$ in $\mathcal{D}$. 
The finite element solution was computed for the same values and the total relative error, given by,
\[ \text{total error} = \sum_{k=1}^{100}\frac{||u^{fe}(\lambda_k) - u_N(\lambda_k)||_{H^1(\Omega)}}{||u^{fe}(\lambda_k)||_{H^1(\Omega)}} \] 
was used as a measure of the accuracy of the solution. 
The computational running time for all four parameter selection methods was also recorded for each value of $N$.
The total relative error and the timing for all the methods (averaged over 10 runs for the greedy and gradient algorithms due to the random selection of the first parameter) are given in Figures~\ref{fig:relerror} and \ref{fig:timing}. 

We see that for all basis sizes, the solutions generated using the Metropolis algorithm have a lower relative error than those generated by the other three algorithms, especially for small basis sizes. The error is almost a full order of magnitude lower than the greedy algorithm for small basis sizes, but was not much lower than the error resulting from use of the gradient algorithm.  
The basis generated by the logarithmic spacing interestingly also mostly outperformed, in terms of relative error, the greedy and gradient algorithms for small basis sizes, confirming the numerical results in~\cite{veroy_thesis}.
We note that the error for the the logarithmic spacing algorithm is not monotonically decreasing,  and posit that this has to do with the specific application, and the tendency of the wavelengths providing the most accurate solutions being clustered together. 

The rapid decline in the relative error for all three methods is expected, as it has been shown numerically that the reduced basis approximation converges rapidly to the finite element solution~\cite{machiels}.
This convergence has been proven to be exponential if the parameter values in $S_N$ are chosen to be logarithmically (quasi-) uniformly~\cite{maday02,prudhomme}.
The choice of basis generating values is fairly uniform for all four methods so the apparent exponential convergence of the relative error in our case is expected. 

Finally, we see that the gradient algorithm has a significantly shorter running time than either of the other two proposed algorithms, discounting the logarithmic spacing that clearly has the fastest running time (practically zero) since no optimality conditions need to be met in this case.
Comparatively, the MCMC algorithm has only a slightly shorter running time than the greedy algorithm for larger basis sizes since the greedy algorithm increases at a slightly faster rate.

First we note that this graph shows the running times for the algorithms given a fixed basis size. If we did not know the optimal or desired basis size up front, we note that the running time for the greedy and gradient algorithms to find a new basis element would be much faster since
they store the previous elements and add one wavelength at a time. This is
in contrast to the Metropolis algorithm which finds an optimal generating
set from scratch for each new N. Although the running time for the Metropolis and greedy algorithms appears to increase linearly with $N$, the gradient algorithm does not monotonically increase. This is unexpected as we assume the runtime to increase with the base size. 

However this behavior can be explained by the fixed base size used in this example and our implementation of the gradient method. As mentioned in Section~\ref{sec:grad}, the gradient method uses a coarse mesh $\Lambda$ to find a starting point for its minimization method. This mesh is initialized with a small number of equidistant points, $N_1$. If the base size $N$ is larger or the same as $N_1$ the mesh is extended to avoid clustering starting points $\mu_k^0$. By providing a fixed base size instead of having a dynamic system the coarse mesh $\Lambda$ this set is different for each $N$ and thus different starting points are chosen. This can be avoided by having a fixed set $\Lambda$ independent of $N$ but bears the problem of possible clustering of starting points which can lead to instability of samples that are to close to each other are being used.

\section{Conclusion}
In this paper, we proposed two alternatives to the greedy algorithm, the Metropolis and gradient algorithms, and compared all three to the method of choosing a basis using a logarithmic spacing. We have demonstrated that both the Metropolis and gradient algorithms may be viable alternatives to the traditional greedy approach when generating the set $S_N$ with which to construct a reduced basis for solving one dimensional parametric PDEs. For our specific one-dimensional parameter ($\lambda$) numerical example, the relative error of the solution generated by the resulting reduced basis for the Metropolis algorithm can be as much as an order of magnitude more accurate than that generated by the greedy algorithm. Additionally, we have shown that the running time of the gradient algorithm is significantly smaller than the other two algorithms which makes a strong case for it being the preferred algorithm when the desired size of the basis is not known since the error of the basis it produces is only slightly larger of that generated by the Metropolis algorithm.

 We have also shown that once the dimension of the reduced basis method is somewhat determined, the Metropolis algorithm works really well and is very promising in finding a global minimum outperforming both greedy and gradient algorithms in our example. 
We do note, however, that the Metropolis algorithm has a very different construction than both the greedy and gradient algorithms, and so a comparison between them may not be completely fair. Empirically, however, using the measures that would matter to a practitioner (accuracy and time), we have demonstrated that the Metropolis algorithm is a viable alternative to the greedy algorithm in 1-D applications such as the one presented here.

We have also demonstrated the efficacy of the proposed algorithms using the forward problem of hyDOT which is very insensitive to both the parameters $D$ and $\mu_a$ as a function of $\lambda$. 
Based on our results, we suspect that the parameter space reduction given by the reduced basis method may help to reduce the complexity of the hyDOT inverse (reconstruction) problem, which is a severely ill-posed problem and thus, may have applications to other inverse problems as well. 
\section{Acknowledgements}
%
Funding: The work of P. Gralla was supported by Deutsche Forschungsgemeinschaft (DFG, German Research Foundation) for subproject C2 ÔOberflŠchenoptimierungÕ within in the SFB 747 (Collaborative Research Center) ÔMikrokaltumformen Ð Prozesse, Charakterisierung, OptimierungÕ. 

\section{References} 
\bibliographystyle{model1b-num-names} 
\bibliography{references}
\end{document}